\documentclass[a4paper,12pt]{article}
\usepackage{hyperref}
\usepackage{stmaryrd}
\usepackage{amsfonts}
\usepackage{bbm}
\usepackage{latexsym,amssymb,amsmath,amscd,amsthm}  
\usepackage[utf8]{inputenc}
\usepackage[T1]{fontenc}
\usepackage{lmodern}

\newtheorem{thm}{Theorem}[section]

\newtheorem{cor}[thm]{Corollary}
\newtheorem{pro}[thm]{Proposition}
\newtheorem{defi}[thm]{Definition}
\theoremstyle{definition}
\newtheorem{ex}[thm]{Example}
\newtheorem{rmk}[thm]{Remark}

\usepackage{color}

\newcommand\strong[1]{\textbf{#1}}  

\newcommand{\be }{\begin{equation}}
\newcommand{\ee }{\end{equation}}

\providecommand\qed{\ensuremath{\Box}}
\renewenvironment{proof}[1][Proof]{\par\medskip\noindent\textbf{#1.} }{\hfill\qed\par\medskip}
\newcommand{\pf}{\begin{proof}} 
\newcommand\Qed{\end{proof}}

\newcommand\ud{{\mathrm{d}}}  
\newcommand\Dsp{{\mathrm{D}}}  
\newcommand\F{{\mathrm{F}}}  

\newcommand{\Real}{{\mathbb R}}
\newcommand{\Comp}{{\mathbb C}}
\newcommand\bi{{\mathsf{i}}}  

\newcommand{\huaG}{{\mathcal{G}}}  
\newcommand\TV{T^V}  

\newcommand\T{{\mathbb{T}}} 

\newcommand{\frkg}{{\mathfrak g}}

\newcommand{\frkp}{{\mathfrak p}}

\newcommand{\Courant}[1]{\left\llbracket  #1\right\rrbracket }

\DeclareMathOperator\im{im}  

\newcommand{\br}[1]{   [ \cdot,    \cdot  ]   }

\DeclareMathOperator\spans{span}

\DeclareMathOperator{\Ad}{Ad}

\newcommand\SO{{\mathrm{SO}}}  
\DeclareMathOperator\Euc{Euc}  

\DeclareMathOperator{\ad}{ad}

\newcommand\email[1]{\href{mailto:#1}{#1}}
\newcommand\eprint[1]{\href{http://arXiv.org/abs/#1}{arXiv: #1}}

\begin{document}
\title{Dirac and Generalized Complex Structures on Cartan Geometries}
\author{Honglei Lang and Xiaomeng Xu \\
   Department of Mathematics and LMAM, Peking University,\\
   Beijing 100871, China\\
  \email{hllang@pku.edu.cn},  \email{xuxiaomeng@pku.edu.cn}
}
\date{}
 \maketitle

\begin{abstract}
We introduce linear Dirac and generalized complex structures on Cartan geometries and give criteria for Dirac subalgebras of $\frkg\ltimes\frkg^*$ representing Dirac structures on a Cartan geometry. We prove that there is a bijection between the linear generalized structures on a torsion free Cartan geometry and the equivariant generalized structures on its corresponding homogeneous space.
\end{abstract}

\section{Introduction}
Generalized complex geometry was introduced in \cite{Hitchin} as a common generalization of complex and symplectic structures and has found several applications in physics and mathematics.
Using the identification between the Lie algebra $\frkg$ of a Lie group $G$ and left invariant vectors field on $TG$, equivariant Dirac structures and generalized complex structures on homogeneous spaces have been classified in \cite{homogeneous} by the Lie algebra information of the homogeneous spaces. This gives a description of equivariant Dirac structures in terms of algebra data.

Cartan geometry is a common generalization of Riemannian geometry, conformal geometry and Klein geometry, which generalizes the linear tangent spaces of the former to the more general homogeneous spaces of the latter. It has been a unifying framework and a powerful tool for the study of conformal, projective, CR and related contact structures in differential geometry. Roughly speaking, a Cartan geometry is a curved analogon of a homogeneous space, twisted by a curvature $\kappa$ (see \cite{3}, \cite{weylstructure}). Every Cartan geometry has a corresponding flat model (homogeneous space $G/P$) and a Cartan connection which plays the role of the Cartan $1$-form in the homogeneous space case. So just like the construction of equivariant Dirac structures in homogeneous spaces by using the Cartan $1$-form, it is natural to introduce certain Dirac structures in a Cartan geometry, which we will call linear, by using the Cartan connection. As a $curved$ homogeneous space, Cartan geometry inherits the generalized structure on its flat model (homogeneous space) in some sense. This gives a description of such Dirac structures in terms of Lie algebra data and allows us to provide some new examples of generalized complex structures.

The paper is organized as follows. In Section~\ref{s:2}, we recall some basic notions of Dirac and generalized complex structures. In Section~\ref{s:3}, we give the basic notions of Cartan geometry including Cartan connection and Cartan curvature. In Section~\ref{s:4}, we introduce linear generalized structures on Cartan geometry and give criteria for a Dirac subalgebra of $\frkg\ltimes\frkg^*$ to represent a Dirac structure on a Cartan geometry. We prove that there is a bijection between the linear generalized structures on a torsion free Cartan geometry and the equivariant generalized structures on its corresponding homogeneous space.

\section{Dirac and Generalized Complex Structures}\label{s:2}
In this section, we introduce the basic definitions of Dirac and Generalized complex structures, see \cite{GualtieriGeneralizedComplex} for more details. For a manifold $M$, there is a natural bilinear form on $\T:=TM\oplus T^*M$, given by pairing
$\langle X+\xi,Y+\eta\rangle_{TM\oplus T^*M}=\eta(X)+\xi(Y)$ for $X,Y\in\Gamma(TM)$ and $\xi,\eta\in\Gamma(T^*M)$. Furthermore, $\T$ is equipped with the Courant bracket defined by
\be
  \Courant{X+\xi,Y+\eta}=[X,Y]+L_X\eta-L_Y\xi-\tfrac12d(\eta(X)+\xi(Y)).
\ee
The bracket $\Courant{\cdot,\cdot}$ and the bilinear form $\langle\cdot,\cdot\rangle$ extend $\Comp$-bilinearly to $\T^\Comp=\T\otimes \Comp$.
\begin{defi}
Given a manifold $M$, a real (complex) almost Dirac structure on $M$ is a maximal isotropic subbundle $\mathbb{D}$ of $\T$ ($\T^\Comp$). A real (complex) almost Dirac structure is called a real (complex) Dirac structure if it is integrable with respect to the Courant bracket.
\end{defi}

\begin{defi}[\cite{GualtieriGeneralizedComplex}]
A generalized almost complex structure on $M$ is a map $J:\T\longrightarrow\T$ such that $J$ is orthogonal with respect to the bilinear form $\langle\cdot,\cdot\rangle$ and $J^2=-1$. A generalized almost complex structure $J$ is called a generalized complex structure if the $\bi$-eigenbundle of $J$ in $\T^\Comp$ is integrable with respect to the Courant bracket.
\end{defi}

It is easy to see that the $\bi$-eigenbundle $\mathbb{D}$ of a generalized complex structure $J$ is actually a maximal isotropic subbundle of $\T^\Comp$, so it is a complex Dirac structure and the following theorem shows that the study of generalized structures lies in the framework of Dirac structures.

\begin{thm}[\cite{GualtieriGeneralizedComplex}]\label{complex}
A generalized (almost) complex structure $J$ is equivalent to a complex (almost) Dirac structure $\mathbb{D}$ such that $\mathbb{D}\cap \overline{\mathbb{D}}=0$.
\end{thm}
It is well-known that for a Dirac structure $\mathbb{D}$ on $TM\oplus T^*M$, the real index $r$ of $\mathbb{D}$ satisfies:
\[r:=\dim_\Comp\mathbb{D}\cap \mathbb{D}^*\equiv m \pmod{2},\]
where $m$ is the dimension of the manifold $M$. So generalized complex structures only exist on even-dimensional manifolds and are actually kinds of minimal real index Dirac structures. For odd-dimensional manifolds, we could also study the minimal real index Dirac structures, i.e.\ $\mathbb{D}\cap \overline{\mathbb{D}}$ is a trivial (under some geometric consideration) subbundle of $TM\oplus T^*M$. They have been researched in some slightly different form under the name generalized contact geometry.

\section{Cartan geometry}\label{s:3}
Roughly speaking, a Cartan geometry is the curved analogon of a
homogeneous space twisted by a curvature $\kappa$.

At first, let us recall the basic definitions in Cartan geometry.
\begin{defi}\label{d:3.1}
Let $G$ be a Lie group with
Lie algebra $\frkg$, a Cartan geometry $(\pi:\huaG\longrightarrow M,\omega,G/P)$ of type $(G,P)$ is
given by a principal bundle $\huaG\longrightarrow{M}$ with structure group $P$ equipped with a
$\frkg$-valued $1$-form (Cartan connection) satisfying the following conditions:
\begin{enumerate}
\item\label{cart:p1}
the map $\omega_p:T_p\huaG\longrightarrow{\frkg}$ is a linear isomorphism for every
$p\in{\huaG}$;

\item\label{cart:p2} ${R_a}^*\,\omega=\Ad(a^{-1})\circ\omega$ $\forall{a}\in{P}$, where $R_a$ denotes the right action of an element $a$ in the structure group;

\item\label{cart:p3} $\omega(\zeta_{A})=A,$ $\forall{A}\in{\frkp}$ where $\frkp$ is the Lie algebra of $P$ and $\zeta_{A}$ is fundamental vector field on $\huaG$ generalized to $A\in\frkg$.
\end{enumerate}

\end{defi}

The \strong{curvature} of a Cartan connection is defined to be a $\frkg$-valued $2$-form $\kappa\in\Gamma(\wedge^2T\huaG\otimes\frkg)$
\be\label{curvature}
\kappa:=\ud\omega+\tfrac12[\omega,\omega].
\ee

Given a homogeneous space $(G,P)$, with the left invariant Cartan $1$-form, we see that they satisfy the conditions of Definition~\ref{d:3.1} and the curvature $\kappa=\ud\omega+\tfrac12[\omega,\omega]=0$ (Maurer--Cartan equation). So homogeneous spaces are flat Cartan geometries, and for a Cartan geometry $(\pi:\huaG\longrightarrow M,\omega,G/P)$, we see the corresponding homogeneous space $(G,P)$ as its flat model.

One simple property is that the curvature $\kappa$ is equivariant (because $\omega$ is equivariant). Hence, we can view $\kappa$ as a $2$-form on $M$ with values in the vector bundle
$\huaG\times_P\frkg$ (associated to the adjoint representation of $G$). Moreover, $\kappa$ is a horizontal form (see \cite{Pontryagin class}), i.e.\ for $\xi\in{\TV P}$ ($\TV P$ is the vertical distribution of principal bundle $(\huaG,P)$), we have $\kappa(\xi,\cdot)=0$. This is to say, we can
view $\kappa$ as a function valued in $\wedge^2(\frkg/\frkp)^*\otimes \frkg$, i.e.
$\kappa\in{C^\infty(\huaG,\wedge^2(\frkg/\frkp)^*\otimes \frkg)}$.

We say that a Cartan geometry is \strong{torsion free} if its curvature function takes values in the subalgebra $\frkp\subset \frkg$ only.

As a particularly interesting case let us consider \strong{parabolic geometry}, i.e.\ a Cartan geometry where $\frkg$ is semisimple and $\frkp$ is a parabolic subalgebra of $\frkg$. We know that (see \cite{Pontryagin class}) for a semi-simple Lie algebra $\frkg$, a parabolic subalgebra $\frkp$ is
equivalent to having an $|l|$-grading of $\frkg$:
\[\frkg=\frkg_{-l}+\cdots+\frkg_{-1}+\frkg_0+\frkg_{1}+\cdots+\frkg_l \]
with $[\frkg_i,\frkg_j]\subset\frkg_{i+j}$, $\frkg_-=\frkg_{-l}+\cdot\cdot\cdot+\frkg_{-1}$ is generated by $\frkg_{-1}$ and we have the following facts about this grading:
\begin{enumerate}
\item\label{gr:pa}
There exists a grading element $E\in\frkg_0$ such that $[E,A]=iA$ if and only if $A\in\frkg_i$;

\item\label{gr:pb} The Killing form $B(\cdot,\cdot)$ defines isomorphisms $\frkg_{-i}\cong \frkg^*_i$ and $B(\frkg_i,\frkg_j)=0$ for all $j\neq-i$.
\end{enumerate}

Now, the curvature of a parabolic geometry can be defined as the function
$\kappa:\huaG\longrightarrow \wedge^2\frkg_-^*\otimes \frkg: p\longmapsto \kappa_p$ with
\[
  \kappa_p(A,B)=\ud\omega_p(\omega_p^{-1}(A), \omega_p^{-1}(B))+[A,B], \quad\forall A,B\in \frkg_.
\]

We can decompose the curvature in terms of this grading, i.e.
\[\kappa_i(p):\wedge^2{\frkg_-}\longrightarrow \frkg_i, \ i=-l,\cdots,l.\]
So a parabolic geometry is torsion free if and only if $\kappa_{-i}=0$, for all $i=1,\dots,l$.

\section{Generalized Structures in Cartan Geometry}\label{s:4}
In this section, we introduce linear generalized structures in Cartan geometry. In the case of homogenous spaces, they are the generalized structures researched in \cite{homogeneous}.

Given a Cartan geometry $(\pi:\huaG\longrightarrow M,\omega,G/P)$, because the map $\omega_p:T_p\huaG\longrightarrow{\frkg}$ is a linear isomorphism for every
$p\in{\huaG}$, we can identify $T\huaG$ with $\huaG\times \frkg$ and denote the inverse by $\omega^{-1}$. Now consider the natural morphism $\F$ from $\huaG\times(\frkg\oplus \frkg^*)$ to $T\huaG\oplus T^*\huaG$ defined by
\[
\F(A+\xi)=\omega^{-1}(A)+\omega^*(\xi),~~ \forall A\in\Gamma(\huaG\times\frkg), \xi\in\Gamma(\huaG\times\frkg^*),
\]
then for every subspace $\Dsp$ of $\frkg\oplus \frkg^*$, we have $\F(\Dsp)$ is a subbundle of $T\huaG\oplus T^*\huaG$, i.e.\ for every point $p$ in $\huaG$
\be\label{F}
\F(\Dsp)_p= \{\omega_p^{-1}(A)+\omega^*(\xi) | A+\xi\in \Dsp \}.
\ee

Define a bilinear form on the sections of $\huaG\times(\frkg\oplus\frkg^*)$ by the natural pairing
\[
\langle A+\xi,B+\eta\rangle_{\huaG\times(\frkg\oplus\frkg^*)}=\eta(A)+\xi(B),
\]
for all $A,B\in\Gamma({\huaG\times\frkg})$ and $\xi,\eta\in\Gamma({\huaG\times\frkg^*})$.

Define a semidirect product Lie algebra structure on the vector space $\frkg\oplus\frkg^*$, i.e.
\[  [A+\xi,B+\eta]_{\frkg\ltimes \frkg^*}=[A,B]_\frkg+{ad^*_A}\eta-{ad^*_B}\xi,  \quad\forall A+\xi,B+\eta\in\frkg\oplus \frkg^*,
\]
and denote it by $\frkg\ltimes \frkg^*$.

\begin{defi}
A subalgebra of $\frkg\ltimes \frkg^*$ is called a Dirac subalgebra if it is a maximal isotropic subspace of $\frkg\oplus \frkg^*$.
\end{defi}
We see that a Dirac subalgebra is just the usual notion of a Dirac structure for the Courant algebroid on $\frkg\oplus\frkg^*$ over a point.

Now we give the criteria for a Dirac subalgebra $\Dsp$ representing a Dirac structure on $\huaG$.
Define $K:\wedge^2(\huaG\times (\frkg\oplus\frkg^*))\longrightarrow \huaG\times \frkg$ by
\[
  K(A+\xi,B+\eta)=\kappa(\omega^{-1}(A),\omega^{-1}(B)),
\] for every $A,B\in\Gamma(\huaG\times \frkg)$, $\xi,\eta\in\Gamma(\huaG\times\frkg^*)$.

Moreover we define a $3$-form $\Theta\in\Gamma(\wedge^3\huaG\times (\frkg\oplus\frkg^*))$ such that for every $e_1,e_2,e_3\in\Gamma(\huaG\times (\frkg\oplus\frkg^*))$
\[
  \Theta(e_1,e_2,e_3)=\langle K(e_1,e_2),e_3\rangle_{\huaG\times(\frkg\oplus\frkg^*)}+\langle K(e_2,e_3),e_1\rangle_{\huaG\times(\frkg\oplus\frkg^*)}+\langle K(e_3,e_1),e_2\rangle_{\huaG\times(\frkg\oplus\frkg^*)}.
\]

\begin{thm}\label{thm:Diracstructure}
Let $(\pi:\huaG\longrightarrow M,\omega,G/P)$ be a Cartan geometry, then for every Dirac subalgebra $\Dsp$ of $\frkg\ltimes \frkg^*$, $\F(\Dsp)$ is a Dirac structure of $T\huaG\oplus T^*\huaG$ if and only if $\Theta|_{\F(\Dsp)}=0$, i.e.,
\be\label{criterion}
\langle K(e_1,e_2),e_3\rangle_{\huaG\times(\frkg\oplus\frkg^*)}+\langle K(e_2,e_3),e_1\rangle_{\huaG\times(\frkg\oplus\frkg^*)}+\langle K(e_3,e_1),e_2\rangle_{\huaG\times(\frkg\oplus\frkg^*)}=0
\ee
for all $e_1,e_2,e_3\in \Dsp$.
\end{thm}

\pf $\F(\Dsp)$ is obviously maximal isotropic, we just need to prove that $[\Gamma(\F(\Dsp)),\Gamma(\F(\Dsp))]\subset \Gamma(\F(\Dsp))$, i.e.\ for every $A+\xi$, $B+\eta\in\Dsp$ (see as constant sections of $\huaG\times(\frkg\oplus\frkg^*)$)
\[
\Courant{\F(A+\xi),\F(B+\eta)}=\Courant{\omega^{-1}(A)+\omega^*(\xi),\omega^{-1}(B)+\omega^*(\eta)}\subset \F(\Dsp).
\]
By formula \eqref{curvature}, i.e.\ $\ud\omega+\tfrac12[\omega,\omega]=\kappa$,  we have for every $A,B,C\in \frkg$ and $\xi,\eta\in \frkg^*$ (seen 
 as constant sections of $\huaG\times \frkg$ and $\huaG\times \frkg^*$ respectively),
\[
\kappa(\omega^{-1}(A),\omega^{-1}(B))=\ud\omega(\omega^{-1}(A),\omega^{-1}(B))+\tfrac12[\omega,\omega](\omega^{-1}(A),\omega^{-1}(B)),
\]
After some simple computations, we get
\be\label{five}
\Courant{\omega^{-1}(A),\omega^{-1}(B)}=\omega^{-1}\big([A,B]_\frkg+k(\omega^{-1}(A),\omega^{-1}(B))\big).
\ee
where $[\cdot,\cdot]_\frkg$ denotes the pointwise Lie bracket of $\huaG\times \frkg$.
On the other hand, it is easy to obtain the following formulas:
\begin{align}
  \nonumber\Courant{\omega^{-1}(A),\omega^*(\eta)}(\omega^{-1}(C))
  &= -\omega^*(\eta)(\omega^{-1}([A,C]_\frkg))-\omega^*(\eta)(\omega^{-1}(k(\omega^{-1}(A),\omega^{-1}(C))))  \\
  &= -\eta([A,C]_\frkg)-\langle\eta,\kappa(\omega^{-1}(A),\omega^{-1}(C))\rangle,  \label{six}\\
  \Courant{{\omega}^*(\xi),{\omega}^*(\eta)}&=0.  \label{seven}
\end{align}
By \eqref{five}, \eqref{six} and \eqref{seven}, we get
\begin{align}
  \nonumber&\Courant{\omega^{-1}(A)+\omega^*(\xi),\omega^{-1}(B)+\omega^*(\eta)} \\
  \nonumber=\,&\omega^{-1}([A,B]_\frkg+\kappa(\omega^{-1}(A),\omega^{-1}(B)))+\omega^*(\ad^*_A\eta)  \\
  \nonumber&\quad-\langle\eta,\kappa(\omega^{-1}(A),\cdot)\rangle-\omega^*(\ad^*_B\xi)+\langle\xi,\kappa(\omega^{-1}(B),\cdot)\rangle  \\
  \nonumber=\,&\F([A+\xi,B+\eta]_{\frkg\ltimes\frkg^*})+\omega^{-1}(\kappa(\omega^{-1}(A),\omega^{-1}(B)))  \\
  &\quad-\langle\eta,\kappa(\omega^{-1}(A),\cdot)\rangle+\langle\xi,\kappa(\omega^{-1}(B),\cdot)\rangle,\label{eight}
\end{align}
where $\langle\eta,\kappa(\omega^{-1}(A),\cdot)\rangle\in\Gamma(T^*\huaG)$ is defined by
\[
  \langle\eta,\kappa(\omega^{-1}(A),\cdot)\rangle(X)=\langle\eta,\kappa(\omega^{-1}(A),X)\rangle,
    \quad\text{for all } X\in\Gamma(T\huaG).
\]

Now assume $A+\xi,B+\eta\in \Dsp$, then $\Courant{A+\xi,B+\xi}\in\F(\Dsp)$ if and only if
\[f(A+\xi,B+\eta):=\omega^{-1}(\kappa(\omega^{-1}(A),\omega^{-1}(B)))-\langle\eta,\kappa(\omega^{-1}(A),\cdot)\rangle+\langle\xi,\kappa(\omega^{-1}(B),\cdot)\rangle \in \F(\Dsp).
\]

It is equivalent to $\langle f(A+\xi,B+\eta),e\rangle_{T\huaG\oplus T^*\huaG}=0$ for every $e\in\Gamma(\F(\Dsp))$.  But this is equivalent to the criterion \eqref{criterion}
.

On the other hand, we have
\[ \Courant{f_1e_1,f_2e_2}=f_1f_2\Courant{e_1,e_2}+f_1(\rho(e_1)f_2)e_2-f_2(\rho(e_2)f_1)e_1
\]
for every $f_1,f_2\in C^\infty(\huaG)$ and $e_1,e_2\in\Gamma( F(D))$.

That is to say that $[\omega^{-1}(A)+\omega^*(\xi),\omega^{-1}(B)+\omega^*(\eta)]\in F(D)$ for every (not necessarily constant) $A+\xi,B+\eta\in\Gamma(\huaG\times(\frkg\oplus\frkg^*))$. This completes the proof.
\Qed

\begin{cor}\label{construct}
Let $(\pi:\huaG\longrightarrow M,\omega,G/P)$ be a torsion free Cartan geometry, then for every Dirac Lie algebra $\Dsp$ of $\frkg\ltimes \frkg^*$ containing $\frkp$, we have $\F(\Dsp)$ is a Dirac structure of $T\huaG\oplus T^*\huaG$.
\end{cor}
\pf Because of torsion freeness, we have $K(e_1,e_2)\in\Gamma(\huaG\times\frkp)$, for all $e_1,e_2\in\Dsp$. Now $\Dsp$ is isotropic and $\im(\kappa_p)\subset \frkp\subset \Dsp$ for every point $p$ on $P$, so $\langle e_1,\frkp\rangle=0$, i.e.\ $\langle e_1,\kappa(e_2,\cdot)\rangle=0,$ for every $e_1,e_2\in\Dsp$. So by Theorem 4.2, the proof is complete.
\Qed

\begin{pro}\label{invariant}
For every Dirac subalgebra $\Dsp$ of $\frkg\ltimes \frkg^*$ containing $\frkp$, $\F(\Dsp)$ is a $P$-invariant subbundle of $T\huaG\oplus T^*\huaG$.
\end{pro}
\pf The invariance of $\F(\Dsp)$ can be obtained from the $P$-equivariance of the Cartan connection, i.e.\ ${R_a}^*\omega=\Ad(a^{-1})\circ\omega$, $\forall{a}\in{P}$, and the fact that $[\frkp,\Dsp]\subset \Dsp$.\Qed

Now let $(\pi:\huaG\longrightarrow M,\omega,G/P)$ be a Cartan geometry, $\Dsp$ be a Dirac subalgebra of $\frkg\ltimes \frkg^*$ containing $\frkp$, and $\F(\Dsp)$ be a Dirac structure on $T\huaG\oplus T^*\huaG$.  By Proposition~\ref{invariant}, we have $\F(\Dsp)$ is a $P$-invariant subbundle.

In order to obtain Dirac structures on $M=\huaG/P$, we use the natural projection $\pi$ to push forward the Dirac structure from $\huaG$ to $M$, i.e.\ define a subbundle $\pi_*(\F(\Dsp))$ of $TM\oplus T^*M$ by
\be\label{pushDirac}
\pi_*(\F(\Dsp)):= \{X+\xi\in TM\oplus T^*M|\widetilde{X}+\pi^*(\xi)\in \F(\Dsp) \},
\ee
where $\widetilde{X}\in\Gamma(TP)$ denotes any $P$-invariant lift of $X\in\Gamma(TM)$.

\begin{thm}\label{thm:Dirac}
Let $(\pi:\huaG\longrightarrow M,\omega,G/P)$ be a Cartan geometry, then for every Dirac Lie algebra $D$ of $\frkg\ltimes \frkg^*$ containing $\frkp$ such that $\F(\Dsp)$ is a Dirac structure, we have $\pi_*(\F(\Dsp))$ is a Dirac structure of $TM\oplus T^*M$.
\end{thm}
\pf At first, we have $\langle e_1, e_2\rangle_{TM\oplus T^*M}=\langle \widetilde{e_1}, \widetilde{e_2}\rangle_{T\huaG\oplus T^*\huaG}$, where $e_1,e_2\in\Gamma(\pi_*\F(\Dsp))$ and $\widetilde{e_1},\widetilde{e_2}\in\Gamma(\F(\Dsp))$ such that $\pi_*(\widetilde{e_i})=e_i$, $i=1,2$. So $\pi_*(\F(\Dsp))$ is obviously isotropic.  Suppose $X+\xi\in\Gamma(TM\oplus T^*M)$ such that
$\langle X+\xi,e\rangle=0$ for every $e\in\Gamma(\pi_*(\F(\Dsp)))$.  Because $\TV \huaG\subset \F(\Dsp)$, where $\TV \huaG$ is the vertical distribution, we get $\langle\widetilde{X}+\pi^*(\xi),e\rangle=0$ for every $P$-invariant lift $\widetilde{X}$ of $X$ and every $e\in\Gamma(\F(\Dsp))$. Because $\F(\Dsp)$ is maximal
, we get $\widetilde{X}+\pi^*(\xi)\in\Gamma(\F(\Dsp))$. By definition of $\F(\Dsp)$, $X+\xi\in\Gamma(\pi_*(\F(\Dsp)))$.

To see the integrability, let $X+\xi$, $Y+\eta$ be sections of $\pi_*(\F(\Dsp))\subset TM\oplus T^*M$ and $\widetilde{X}+\pi^*(\xi)$, $\widetilde{Y}+\pi^*(\eta)$ be the corresponding sections of $\F(\Dsp)\subset TP\oplus T^*P$.  We have
\[
  \Courant{\widetilde{X}+\pi^*(\xi),\widetilde{Y}+\pi^*(\eta)}=[\widetilde{X},\widetilde{Y}]
   +L_{\widetilde{X}}\pi^*(\eta) -L_{\widetilde{Y}}\pi^*(\xi).
\]
But $L_{\widetilde{X}}\pi^*(\eta)-\widetilde{L_X\eta}\in \TV P$ and $[\widetilde{X},\widetilde{Y}]-\widetilde{[X,Y]}\in \TV P$ where $\widetilde{L_X\xi}$ and $\widetilde{[X,Y]}$ denote any $P$-invariant lift of $L_X\xi$ and $[X,Y]$ respectively.  So $\widetilde{[X,Y]}+\pi^*(L_X\eta)-\pi^*(L_Y\xi)\in\Gamma(\F(\Dsp))$. Thus by definition $\Courant{X+\xi,Y+\eta}\in\Gamma(\pi_*(\F(\Dsp)))$. This completes the proof.
\Qed

Summarizing Theorem~\ref{thm:Diracstructure} and Theorem~\ref{thm:Dirac}, we get the main theorem of this paper:
\begin{thm}\label{main}
Let $(\pi:\huaG\longrightarrow M,\omega,G/P)$ be a Cartan geometry with Cartan curvatrue $\kappa$, then $\pi_*(\F(\Dsp))$ gives a Dirac structure over $M$ if $\frkp\subset\Dsp$ and $\langle K(e_1,e_2),e_3\rangle+\langle K(e_2,e_3),e_1\rangle+\langle K(e_3,e_1),e_2\rangle=0$ for all $e_1,e_2,e_3\in\Dsp$.
\end{thm}
\begin{rmk}
By considering the complexifications of $T\huaG\oplus T^*\huaG$, $\frkg$ and $\frkp$, the Theorem~\ref{thm:Dirac} also holds in the complex case.
\end{rmk}

\begin{cor}\label{poisson}
Let $(\pi:\huaG\longrightarrow M,\omega,G/P)$ be a torsion-free Cartan geometry, then for every Dirac subalgebra $D$ of $\frkg\ltimes \frkg^*$ containing $\frkp$, we have $\pi_*(\F(\Dsp))$ is a Dirac structure of $TM\oplus T^*M$.
\end{cor}

It is well-known that Poisson structures are special cases Dirac structures.  Following the discussions in \cite{liuwx}, we have:
\begin{cor}\label{Poisson}
Let $(\pi:\huaG\longrightarrow M,\omega,G/P)$ be a Cartan geometry with Cartan curvatrue $\kappa$, then $\pi_*(\F(\Dsp))$ gives a Poisson structure on $M$ if $\Dsp\cap\frkg=\frkp$, and $\langle K(e_1,e_2),e_3\rangle+\langle K(e_2,e_3),e_1\rangle+\langle K(e_3,e_1),e_2\rangle=0$ for all $e_1,e_2,e_2\in\Dsp.$
\end{cor}
\begin{ex}
Riemannian geometry can be defined as a torsion-free Cartan geometry modeled on an Euclidean space $G/P$ with $G= \Euc_n\Real$, the group of Euclidean motions, and $P=\SO(n)$, see \cite{Riemannian} for more details. In order to obtain a Dirac structure on $M$ by Corollary \ref{poisson}, we just need to find a Dirac subalgebra $\Dsp$ of $\frkg\ltimes \frkg^*$ containing $\frkp$, where $\frkg$ and $\frkp$ are Lie algebras of Lie groups $G$ and $P$ respectively. By Corollary \ref{Poisson}, a Dirac subalgebra $\Dsp$ of $\frkg\ltimes \frkg^*$ gives a Poisson structure on a Riemannian manifold $M$ if $\Dsp\cap\frkg=\frkp$. Now, let $M$ be a $2$-dimensional Riemannian manifold, then its flat model is just $\Euc_2\Real/\SO(2)$. The Lie algebra $\frkg$ of $\Euc_2\Real$ is generated by $e_1,e_2,e_3$, and the Lie bracket is given by $[e_1,e_2]=-e_3$, $[e_1,e_3]=e_2$, $[e_2,e_3]=0$. The Lie subalgebra $\frkp$ of $\frkg$ is generated by $e_1$. Then we choose a subspace $\Dsp$ of $\frkg\ltimes \frkg^*$ by $\Dsp=\spans\{e_1,e_2-e_3^*,e_3+e_2^*\}$. It is easy to see that $\Dsp$ is a Dirac subalgebra and $\Dsp\cap\frkg=\frkp$, so $\Dsp$ gives a Poisson structure on the Riemannian manifold $M$.
\end{ex}

From Theorem \ref{complex} and Corollary \ref{poisson}, we have:
\begin{cor}
Let $(\pi:\huaG\longrightarrow M,\omega,G/P)$ be a torsion free Cartan geometry. If
$\Dsp\subset \frkg_\Comp\ltimes{\frkg_\Comp}^*$ is a Dirac subalgebra containing $\frkp_\Comp$ and
$\F(\Dsp)$ is a Dirac structure on $\huaG$, then $\Dsp$ induces a
generalized complex structure on $M$ if and only if $\Dsp\cap \overline{D}=\frkp_\Comp$.
\end{cor}

\begin{defi}
Let $(\pi:\huaG\longrightarrow M,\omega,G/P)$ be a Cartan geometry, a Dirac structure $L$ of $TM\oplus T^*M$ is called a \strong{linear Dirac structure} if there exists a Dirac subalgebra $\Dsp$ of $\frkg\rtimes \frkg^*$ such that $L=\pi_*(\F(\Dsp))$. A generalized complex structure $J$ on $M$ is called a \strong{linear generalized complex structure} if the $\bi$-eigenbundle of $J$ is a linear Dirac structure.
\end{defi}

As a special case of Theorem \ref{main}, we study the Dirac structures on a homogeneous space $G/P$, where $P$ is a closed connected subgroup. Here the Cartan connection $\omega$ is just the Cartan $1$-form and the curvature $\kappa=0$.  We identify left-invariant vector fields with $\frkg$ by $\omega$. So we see that for every Dirac Lie algebra $\Dsp$ of $\frkg\ltimes \frkg^*$, $\F(\Dsp)$ is the corresponding left invariant Dirac structure on $TG\oplus T^*G$. By Corollary \ref{poisson}, every Dirac Lie algebra $\Dsp$ containing $\frkp$ gives a Dirac structure $\pi_*(\F(\Dsp))$ of $TM\oplus T^*M$ and it is easy to see that this Dirac structure is $G$-invariant. On the other hand, we have.
\begin{thm}[\cite{homogeneous}]\label{D}
There is a bijection between the $G$-invariant Dirac structures on $G/P$ and the set of Dirac subalgebras $\Dsp$ of $\frkg\ltimes \frkg^*$ containing $\frkp$.
\end{thm}
So in the case of a homogeneous space, the Dirac structures and generalized complex structures studied in \cite{homogeneous} are just the linear Dirac structures we get by Theorem \ref{main}.
\par\bigskip

From the discussion above, we see that linear generalized structures are determined by data in the Lie algebra.  At the same time, the $G$-invariant generalized structures on the homogeneous space $G/P$ are determined by the same algebra data (see Corollary \ref{algebracondition} and \cite{homogeneous}). That is to say, there is a bijection between the linear generalized structures on a torsion free Cartan geometry $(\pi:\huaG\longrightarrow M,\omega,G/P)$ and the $G$-invariant generalized structures on its corresponding homogeneous space $G/P$. Every invariant generalized structure on $G/P$ has a corresponding ``curved'' generalized structure on $\huaG/P$, so as a $curved$ homogeneous space, every Cartan geometry inherits the generalized structures of its flat model (homogeneous space) in the above sense.
\par\bigskip

Now, every $G$-invariant generalized structure in torsion free Cartan geometry can be decoded by the corresponding Dirac Lie algebra of $\frkg\ltimes \frkg^*$, i.e.\ linear Dirac structures. On the other hand, every maximal isotropic subspace of $\frkg\ltimes \frkg^*$ is of the form
\[L(E,\varepsilon):= \{X+\xi\in \frkg\oplus \frkg^*|X\in E \ and \ i_X\varepsilon=\xi_{|E} \}
\]
where $E$ is a subspace of $\frkg$ and $\varepsilon\in \wedge^2\frkg^*$.
\begin{pro}
$L(E,\varepsilon)$ is integrable with respect to the semidirect Lie algebra $\frkg\ltimes \frkg^*$ if and only if
\begin{enumerate}
\item\label{int:p1} $E\subset \frkg$ is a subalgebra,
\item\label{int:p2} $\ud_E\varepsilon=0$,
\end{enumerate}
where $\ud_E$ is given by
\[\ud_E\varepsilon(X,Y,Z)=\varepsilon(X,[Y,Z])+\varepsilon(Y,[X,Z])+\varepsilon(Z,[X,Y]),\]
for all $X,Y,Z\in\frkg$.
\end{pro}

\begin{pro}
Let $L(E,\varepsilon)$ be a Dirac subalgebra of $\frkg\ltimes \frkg^*$, then $\frkp\subset L(E,\varepsilon)$ if and only if $\varepsilon^\sharp$ vanishes on $\frkp$.
\end{pro}
By the propositions above and Theorem \ref{main}, we have:
\begin{cor}\label{algebracondition}
The linear Dirac structures on a Cartan geometry $(\pi:\huaG\longrightarrow M,\omega,G/P)$ are parameterized by pairs $(E,\varepsilon)$, where $E$ is a subalgebra of $\frkg$ containing $\frkp$, $\ud_E\varepsilon=0$, $\varepsilon^\sharp$ vanishes on $\frkp$, and $\langle K(e_1,e_2),e_3\rangle+\langle K(e_2,e_3),e_1\rangle+\langle K(e_3,e_1),e_2\rangle=0$ for every $e_1,e_2,e_3\in L(E,\varepsilon)$.
\end{cor}

Following the discussion in \cite{homogeneous}, we get:
\begin{cor}\label{algebracondition2}
Let $(\pi:\huaG\longrightarrow M,\omega,G/P)$ be a Cartan geometry. There is a bijection between linear generalized complex structures on $M$ and pairs $(E,\varepsilon)$ where $E$ is a subalgebra of $\frkg_\Comp$ and $\varepsilon\in\wedge^2E^*$ with
\begin{enumerate}
\item\label{p1} $\frkp_c\subset E$,

\item\label{p2} $E+\overline{E}=\frkg_\Comp$,

\item\label{p3} $\ud_E\varepsilon=0$,

\item\label{p4} $\varepsilon^\sharp(\frkp)=0$,

\item\label{p5} For all $X\in E\cap \overline{E}$, if $\varepsilon(X,Y)-\overline{\varepsilon(\overline{X},\overline{Y})}=0$ for every $Y\in E\cap\overline{E}$, then $X\in\frkp_\Comp$.

\item\label{p6} $\langle K(e_1,e_2),e_3\rangle+\langle K(e_2,e_3),e_1\rangle+\langle K(e_3,e_1),e_2\rangle=0$, for all $e_1,e_2,e_3\in L(E,\varepsilon)$.
\end{enumerate}
\end{cor}
 From Corollary \ref{algebracondition} or \ref{algebracondition2}, we see that for a general Cartan geometry, the criteria for a Dirac subalgebra of $\frkg\ltimes\frkg^*$ to represent a generalized structure are the same as in the homogeneous case except for the extra condition involving the Cartan curvature.  Conversely every Dirac subalgebra of $\frkg\ltimes\frkg^*$ which gives a Dirac structure on a Cartan geometry also gives a $G$-invariant Dirac structure on the corresponding homogeneous space. The homogeneous space case has been studied thoroughly in \cite{homogeneous}, so it would be interesting to analyze some special (torsion free for example) Cartan geometries to obtain new generalized complex structures.

\end{document}